\documentclass[11pt]{amsart}

\usepackage{amsfonts,amsmath,amssymb,amsthm}
\usepackage[T1]{fontenc}

\newcommand{\Cc}{\mathbb{C}}  

\newcommand{\Rr}{\mathbb{R}}
\newcommand{\Nn}{\mathbb{N}}

\newcommand{\defi}[1]{\emph{#1}}
\newcommand{\supp}{\mathop{\mathrm{supp}}\nolimits}

\renewcommand{\epsilon}{\varepsilon}
\renewcommand{\le}{\leqslant}
\renewcommand{\ge}{\geqslant}

{\theoremstyle{plain}
\newtheorem{theorem}{Theorem}    
\newtheorem{lemma}[theorem]{Lemma}        
}
{\theoremstyle{remark}
\newtheorem*{remark*}{Remark}  

\newtheorem{question}{Question}
}

\begin{document}

\title{Jump of Milnor numbers}
\author{Arnaud Bodin}
\address{Laboratoire Paul Painlev\'e, Math\'ematiques,
Universit\'e de Lille I, 
 59655 Villeneuve d'Ascq, France.}
\email{Arnaud.Bodin@math.univ-lille1.fr}
\date{\today}


\begin{abstract}
In this note we study a problem of A'Campo about the minimal non-zero 
difference between the Milnor numbers of a germ of plane curve and one of its deformation.
\end{abstract}

\maketitle

\section{Problem of the jump (A'Campo)}

Let $f_0 : (\Cc^n,0) \longrightarrow (\Cc,0)$ be an analytic germ of isolated singularity.
A \defi{deformation} of $f_0$ is a family $(f_s)_{s\in[0,1]}$ of germs of isolated singularities
such that the coefficients are analytic functions of $s\in[0,1]$.

The \defi{jump of the family $(f_s)$} is 
$$\mu(f_0)-\mu(f_s), \quad 0< s \ll 1,$$
where $\mu$ is the Milnor number at the origin.
This number is well-defined because $\mu(f_0)-\mu(f_s)$ is independent of $s$
if $s$ is sufficiently small, moreover by the upper semi-continuity of $\mu$ this number
is a non-negative integer.

The most famous result about the Milnor number and the topology of the family
is L\^e-Ramanujam's theorem \cite{LR}:
\begin{theorem}
If $n\not=3$ and if $\mu(f_0) = \mu(f_s)$ for all $s\in [0,1]$ then the topological types 
of $f_0^{-1}(0)$ and $f_s^{-1}(0)$ are equal. 
\end{theorem}
In other words, if the jump of the family $(f_s)$ is $0$ then
$f_0^{-1}(0)$ and $f_s^{-1}(0)$ have the same topological type for 
sufficiently small $s$.
Another motivation is that the jump of a family is crucial in the theory of
singularities of polynomial maps at infinity.

The \defi{jump} $\lambda(f_0)$ of $f_0$ is the minimum of the non-zero 
jumps of the $(f_s)$ over all deformations of $f_0$. 
The problem, asked by N.~A'Campo, is to compute $\lambda(f_0)$.
We will only deal with plane curve singularities, that is to say $n=2$.

As a corollary of our study we prove the following:
\begin{theorem}
\label{th:irr}
If $f_0$ is an irreducible germ of plane curve and is Newton non-degenerate
then 
$$\lambda(f_0)=1.$$
\end{theorem}

A closely related question of V.~Arnold \cite{Ar} formulated with our
definitions is to find all singularities with $\lambda(f_0)=1$.
S.~Gusein-Zade \cite{GS} proved that there exist singularities with
$\lambda(f_0)>1$ and as a corollary of a studied of the behaviour of
the Milnor number in a deformation of a desingularization he proved
Theorem \ref{th:irr} for all irreducible plane curves.

This note is organized as follows, in paragraphs \ref{sec:def} to \ref{sec:ex} we define and calculate a
weak form of the jump : the non-degenerate jump. 
In paragraph \ref{sec:irr} we prove Theorem \ref{th:irr} and in paragraph \ref{sec:nonirr}
we give estimations when the germ is not irreducible.
Finally in paragraph \ref{sec:conj} we state some conjectures for the jump of $x^p-y^q$, $p,q\in \Nn$
and end by questions.

\medskip

\emph{Acknowledgments :} I am grateful to Norbert A'Campo for discussions.

\section{Ku\v snirenko's formula}
\label{sec:def}

We firstly recall some definitions (see \cite{K}).
Let $f(x,y) = \sum_{(i,j) \in \Nn^2 } a_{i,j}x^iy^j$ be an analytic germ of plane curve.
Let $\supp(f) = \{ (i,j)\in \Nn^2 \mid a_{i,j} \not=0 \}$
and $\Gamma_+(f)$ be the convex closure of $\bigcup_{(i,j)} ((i,j)+\Rr^2_+)$ where $(i,j)\in \supp(f)\setminus\{(0,0)\}$. The \defi{Newton polygon} $\Gamma(f)$ is the union of the compact faces (called the \defi{slopes}) of  $\Gamma_+(f)$. We often identify a pair $(i,j)\in \Nn^2$ with the monomial $x^iy^j$.
Let $f$ be \defi{convenient} if $\Gamma(f)$ intersects both $x$-axis and $y$-axis.

For a face $\gamma$ of $\Gamma(f)$, let $f_\gamma = \sum_{(i,j) \in \gamma} a_{i,j}x^iy^j$. Then
$f$ is \defi{(Newton) non-degenerate} if for all faces $\gamma$ of $\Gamma(f)$ the system
$$\frac{\partial f_\gamma}{\partial x}(x,y) = 0\ ; \quad \frac{\partial f_\gamma}{\partial y}(x,y) = 0$$
has no solution in $\Cc^*\times \Cc^*$.

For a Newton polygon $\Gamma(f)$, let $S$ be the area bounded by the polygon
 and $a$ (resp. $b$) the length
of the intersection of $\Gamma(f)$ with the axes $x$-axis (resp. $y$-axis).
We set 
$$\nu(f) = 2S-a-b+1.$$
For a convenient germ $f$ the local Milnor number verifies \cite{K} :
\begin{theorem} \ 
\begin{itemize}
\item $\mu(f) \ge \nu(f)$,
\item if $f$ is non-degenerate then $\mu(f)=\nu(f)$.
\end{itemize}
\end{theorem}

\section{Non-degenerate jump for curve singularities}
\label{sec:ndj}

We will consider a weaker problem:
Let $f_0$ be a plane curve singularity  and we suppose that 
$(f_s)$ is a \defi{non-degenerate deformation} that is to say for all 
$s\in]0,1]$, $f_s$ is Newton non-degenerate.
The \defi{non-degenerate jump} $\lambda'(f_0)$ of $f_0$ 
is the minimum of the non-zero jumps over all non-degenerate
deformations of $f_0$.
The new problem is to compute $\lambda'(f_0)$, in this note we
explain how to compute it.

Obviously we have $\lambda(f_0) \le \lambda'(f_0)$ but this inequality can be strict.
For example let $f_0(x,y)=x^4-y^4$, then $\lambda'(f_0)=3$ which is 
obtained for the family $f_s'(x,y)=x^4-y^4+sx^3$. But
$\lambda(f_0) \le 2$, by the degenerate family 
$f_s(x,y)= x^4 - (y^2+sx)^2$ of jump $2$.

\section{Computation of the non-degenerate jump}
\label{sec:comp}

For a convenient $f_0$ there exists a finite set $\mathcal{M}$ of monomials $x^py^q$ lying between the axes 
(in a large sense) and the Newton polygon $\Gamma(f_0)$ (in a strict sense).

\begin{lemma}
If $f_0$ is non-degenerate and convenient then
$$\lambda'(f_0) = \min_{x^py^q\in \mathcal{M}} \big( \mu(f_0) - \mu(f_0+sx^py^q) \big),$$
for a sufficiently small $s\not=0$ (the minimum is over the non-zero values).
\end{lemma}

\begin{proof}
The proof is purely combinatoric and is inspired from \cite{Bo}.
For any polygon $T$ of $\Nn\times \Nn$, we define as for $\nu$ a number $\tau(T) = 2S-a-b$.
Then $\tau$ is additive : let  $T_1, T_2$ be polygons whose vertices are in $\Nn\times \Nn$,
and such that $T_1\cap T_2$ has null area then
$\tau(T_1 \cup T_2) = \tau(T_1) + \tau(T_2)$. By this additivity we can argue on triangles only.
Moreover for a polygon $T$ that do not contain $(0,0)$ we have $\tau(T) \ge 0$.

Now the jump for a non-degenerate family $(f_s)$ corresponds to $\tau(T)$ where 
$T$ is the polygon ``between'' $\Gamma(f_0)$ and $\Gamma(f_s)$ ($0< s \ll 1$).
Minimizing this jump is equivalent to minimizing $\tau(T)$. It is obtained for a polygon $T$
for which all vertices except one are in $\Gamma(f_0)$ and the last vertex is in $\Gamma(f_s)$.
Then it is sufficient to add only one monomial corresponding to the latter vertex
to obtain the required deformation.
\end{proof}

With this method we do not compute $\mu(f_0)$, nor  $\mu(f_0+sx^py^q)$
but directly the difference.

For a degenerate function $f$ we denote by $\tilde f$ a non-degenerate function
such that $f$ and $\tilde f$ have the same Newton polygon: $\Gamma(f_0)=\Gamma(\tilde f_0)$.
The non-degenerate jump for a degenerate function $f_0$ can be computed with
the easy next lemma:
\begin{lemma}
Let $f_0$ be degenerate.
\begin{itemize}
  \item $\lambda'(f_0) =\mu(f_0)-\mu(\tilde f_0)$ if $\mu(f_0)-\mu(\tilde f_0) >0$, 
  \item else $\lambda'(f_0) = \lambda'(\tilde f_0)$.
\end{itemize}
\end{lemma}

\section{An example}
\label{sec:ex}

For a given polynomial $f_0$ it is very fast 
to see who will be the good candidates $x^py^q$ and hence to find $\lambda'(f_0)$
after a very few calculus: we use that $\mu(f_0) - \mu(f_0+sx^py^q) = \tau(T)$ where
 $T$ is the zone between the Newton polygon of $f_0$ and the one of $f_0+sx^py^q$.

For example let $f_0(x,y) = x^4-y^3$.
We draw its Newton polygon (see Figure \ref{fig:ex}). We easily see that the monomials $x^py^q$
that are candidates to minimize $\tau$ for the zone between the Newton polygons are $x^3$ 
(that will give a zone with $\tau(T) = 2$) and $xy^2$ that will give a zone with $\tau(T) = 1$.
In that case the deformation will be 
$f_s(x,y) = x^4-y^3+sxy^2$ and the jump of $f_0$ is $1$.

\begin{center}
\begin{figure}[ht]

\unitlength 1mm

\begin{picture}(50,40)(0,0)

\put(0,0){\vector(1,0){50}}
\put(0,0){\vector(0,1){40}}
\put(50,2){\makebox(0,0){$x$}}
\put(2,40){\makebox(0,0){$y$}}

\multiput(0,0)(10,0){5}{\line(0,1){30}}
\multiput(0,0)(0,10){4}{\line(1,0){40}}

\put(40,0){\circle*{2}}
\put(40,0){\line(-4,3){40}}
\put(0,30){\circle*{2}}

\put(30,0){\circle{2}}
\put(0,30){\line(1,-1){10}}
\put(10,20){\circle{2}}
\put(10,20){\line(3,-2){30}}
\end{picture}

\caption{\label{fig:ex} Example $f_0(x,y)=x^4-y^3$}
\end{figure}
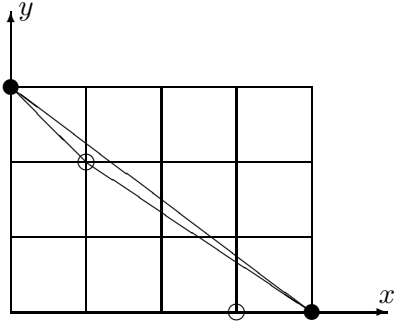
\end{center}

\section{Irreducible case}
\label{sec:irr}

In some cases we are able to give a formula for the computation of the jump.
For example if $f_0(x,y) = x^p-y^q$, with $\gcd(p,q)=d$ then by B\'ezout theorem there exists
a pair $(a,b)$ such that $x^ay^b$ is in $\mathcal{M}$ and such that
the area $T$ corresponding to the deformation $f_s(x,y) = x^p-y^q+sx^ay^b$
is equal to $d/2$.

As an application we prove Theorem \ref{th:irr} cited in the introduction.
\begin{theorem}
If $f$ is irreducible and non-degenerate then
$\lambda(f)=\lambda'(f)=1$.
\end{theorem}

\begin{proof}
We recall some facts from the book of Brieskorn-Kn\"orrer \cite[p.477]{BK}.
For a germ of curve $f$, the number of slopes of a Newton polygon $\Gamma(f)$ is lower or equal to
the number $r$ of irreducible components.

Moreover let $R$ be the number of lattice points that belongs to $\Gamma(f)$ minus $1$.
Then if $f$ is non-degenerate we have $R=r$. The non-degenerate condition
is not explicit in \cite{BK} but it is stated with an equivalent condition
(a face is non-degenerate if and only if the corresponding polynomial $g_i$
of \cite[p.478]{BK} has only simple roots).

Then for an irreducible singular germ $f$, $\Gamma(f)$ has only one slope and $f$ is convenient;
moreover if $f$ is non-degenerate then the extremities of $\Gamma(f)$, say
$x^p$ and $y^q$, verify $\gcd(p,q)=1$.
The non-degenerate jump of $f$ is the same as for $f_0=x^p-y^q$ and is equal to $1$ by B\'ezout theorem.
Then $\lambda'(f)= \lambda'(f_0)=1$, as $0 < \lambda(f) \le \lambda'(f)=1$ it implies
$\lambda(f)=1$.
\end{proof}

\section{Non irreducible case}
\label{sec:nonirr}

More generally if $f$ is convenient, non-degenerate, with one slope, let $x^p$, $y^q$ be the extremities of the Newton polygon of $f$. Then $f$ has the same non-degenerate jump 
as $f_0=x^p-y^q$, we suppose $p\ge q$
and we set $d=\gcd(p,q)$. The formula for $\lambda'(f_0)$ is given by:
\begin{enumerate}
\item If $1 \le d < q \le p$ then $\lambda'(f_0)= d$ which is reached by
a family $f_s(x,y) = x^p-y^q+sx^ay^b$, $a,b$ given by B\'ezout theorem.

\item If $\gcd(p,q)=q$, i.e. $d=q$ then $\lambda'(f_0) = q-1$ which is reached
with $f_s(x,y) = x^p-y^q+sx^{p-1}$.
\end{enumerate}

We will give in paragraph \ref{sec:conj} a conjectural value for $\lambda(x^p-y^q)$.

If there are several slopes with $f$ convenient and non-degenerate then we can estimate $\lambda'(f)$.
Let $f = \prod_{i=1}^k f_i$ be the decomposition of $f$ according to the slopes of
$\Gamma(f)$ (notice that $f_i$ is not necessarily irreducible).
If $f_i$ is a smooth germ then we set (by convention)
$\lambda'(f_i)=1$. In fact $f_i$ is smooth if and only if the corresponding slope $\Gamma_i$ with
extremities $A_i$, $B_i$ verifies $|x_{B_i}-x_{A_i}| = 1$ or $|y_{B_i}-y_{A_i}| = 1$.
Then the following can be proved:
\begin{lemma}
Let $f$ be a convenient non-degenerate germ with several slopes, let $f = \prod_{i=1}^k f_i$
 be the decomposition according to the slopes.
\begin{enumerate}
  \item If all the $f_i$ are smooth then $\lambda'(f)=1$.
  \item If none of the $f_i$ is smooth then
$$ \min_{i=1..k} \lambda'(f_i)  \le \lambda'(f) \le \max_{i=1..k} \lambda'(f_i).$$
  \item In the other cases we have
$$ \min_{i=1..k} \lambda'(f_i)  \le \lambda'(f) \le \max_{i=1..k} {\lambda'(f_i)}+1.$$
\end{enumerate}
\end{lemma}

We give some examples:
\begin{enumerate}
  \item The family $f_s(x,y)=(x+y^4)(x+y^2)(x^2+y)+sy^4$ is of non-degenerate jump $1$.
  \item The family $f_s(x,y)= (x^8-y^6)(x^3-y^2)+sxy^7$ gives $\lambda'(f_0)=2$
with $\lambda'(x^8-y^6)=2$ and $\lambda'(x^3-y^2)=1$.
  \item The family $f_s(x,y)= (x^8-y^6)(x^3-y^2)(x^4-y^4)+sx^5y^7$ verifies $\lambda'(f_0)=2$
while $\lambda'(x^8-y^6)=2$ and $\lambda'(x^3-y^2)=1$ and $\lambda'(x^4-y^4)=3$.
  \item The family $f_s(x,y)= (x+y^3)(x^4+y^4)(x^2+y)+sy^5$ verifies $\lambda'(f_0)=4$
with the smooth germs $x+y^3$, $x^2+y$ and $\lambda'(x^4+y^4) = 3$.
\end{enumerate}

\section{Conjectures for the jump}
\label{sec:conj}

We give a conjectural value for $\lambda(f_0)$ in the case that 
$f_0 = x^p-y^q$ with $p\ge q$.

\begin{enumerate}
\item If $\gcd(p,q)=1$ then $\lambda(f_0)= \lambda'(f_0)=1$.

\item If $p=q$ and $q$ is prime then $\lambda'(f_0)= q-1$, 
with the family $f_s(x,y) = x^q+y^q+sx^{q-1}$. And we conjecture 
that $\lambda(f_0)=q-2$ with the family  $f_s(x,y) =x^q+y^q+s(x+y)^{q-1}$.

\item \label{it:2}
If  $p=kq$ ($k > 1$) and $q$ is prime, then $\lambda'(f_0)= q-1$,
with the family $f_s(x,y) = x^p+y^q+sx^{p-1}$.
It is conjectured that $\lambda(f_0) = \lambda'(f_0)$.

\item \label{it:3}
If $q$ is not prime and  $p=kq$, $k\in \Nn^*$ then let
$q=ab$ with $a\ge 2$ the smallest prime divisor of $q$.
Then $\lambda'(f_0) = q-1 = ab-1$ for the family $f_s(x,y) = x^p-y^q+sx^{p-1}$.
It is conjectured that $\lambda(f_0)= ab-b$, which jump is reached for the
family $f_s(x,y) = x^p - (y^a+sx^{ka-1})^b$.  

\item  \label{it:4}
If $\gcd(p,q)=d$ with $1 < d < q\le p$ then
$\lambda'(f_0) = d$. And it is conjectured that  $\lambda(f_0)=d$ too.

\end{enumerate}

We make a  remark for point (\ref{it:3}), let $g_0(x,y) = x^{p/b}-y^{q/b}= x^{ka}-y^a$.
Then $g_0$ verifies the hypotheses of point (\ref{it:2}) where we have
conjectured $\lambda(f_0)= a-1$ for the deformation
$g_s(x,y) = x^{ka}-y^a-sx^{ka-1}$. Then we calculate 
$g_s(x,y)^b = (x^{ka}-y^a-sx^{ka-1})^b$ which is of course not a reduced polynomial.
We develop and we have an approximation of $g_s(x,y)^b$ if we set
$f_s(x,y) = x^{kab} - (y^a+sx^{ka-1})^b = x^p - (y^a+sx^{ka-1})^b$ with a jump
equal to $ab-b$.

\medskip

Apart from the conjectures above  we ask some questions.
Even if it seems hard to give a formula for the jump, maybe the following
is easier:
\begin{question}
Find an algorithm that computes $\lambda$.
\end{question}

\medskip

Finally the problem of the jump can be seen as a weak form of the problem of adjacency.
For example the list of possible Milnor numbers arising from deformations of $f_0(x,y)=x^4-y^4$ is
$(9,7,6,5,4,3,2,1,0)$. Then the gap between the first term $9=\mu(f_0)$ and the second term is the jump
$\lambda(f_0)=2$.  Then the following question is a generalization of the problem of the jump.
\begin{question}
Give the list of all possible Milnor numbers arising from deformations of a germ.
\end{question}


\end{document}